\documentclass[a4paper]{elsart}

\usepackage{latexsym}
\usepackage{amsmath}
\usepackage{times}

\usepackage{pslatex}
\usepackage{amssymb}
\usepackage{mathrsfs}
\usepackage{cancel}
\usepackage{psfrag, enumerate}

\usepackage{algorithm}
\usepackage{algorithmic}

\usepackage{multicol}

\usepackage[pdftex]{graphicx}



\begin{document}

\begin{frontmatter}
\title{Corotational formulation for 3d solids. An analysis of geometricaly nonlinear foam deformation.}
\author{{\L}ukasz Kaczmarczyk\corauthref{cor}}\ead{Lukasz.Kaczmarczyk@glasgow.ac.uk},
\author{Tomasz Koziara}
and
\author{Chris J. Pearce}
\corauth[cor]{Corresponding author.}
\address{
Department of Civil Engineering, University of Glasgow, Rankine Building
}

\begin{abstract} 

This paper presents theory for the Lagrange co-rotational (CR) formulation of
finite elements in the geometrically nonlinear analysis of 3D structures. In
this paper strains are assumed to be small while the magnitude of rotations
from the reference configuration is not restricted. A new best fit rotator and
consistent spin filter are derived. Lagrange CR formulation is applied with
Hybrid Trefftz Stress elements, although presented methodology can be applied
to arbitrary problem formulation and discretization technique, f.e. finite
volume methods and lattice models, discreet element methods. Efficiency of CR
formulation can be utilized in post-buckling stability analysis, damage and
fracture mechanics, modelling of dynamic fragmentation of bodies made from
quasi-brittle materials, solid fluid interactions and analysis of post-stressed
structures, discreet body dynamics.

\end{abstract}


\end{frontmatter}

\section{Introduction}

Nowadays the Lagrangian kinematics description for Finite Element Method of
geometric nonlinear structures is utilized mainly in two solution techniques,
i.e. Total Lagrangian Formulation and Updated Lagrangian Formulation.  In this
paper we investigate a less common solution technique for geometrically nonlinear
problems, i.e. the Lagrange Co-rotational (CR) formulation.

Key features of the Lagrange CR formulation are efficiency and robustness, i.e.
computation of tangent matrices is is low-priced, the method is easy to parallelize
and can be used with various discretization techniques. Large class of
problems with finite rotations but small strains, e.g. post-buckling
stability analysis, damage and fracture mechanics, modelling of dynamic
fragmentation of bodies made from quasi-brittle materials, solid fluid
interaction, analysis of post-stressed structures or discreet body dynamics
can be solved efficiently with the use of Co-rotational formulation.  

As noted above, the Lagrange CR formulation can be applied together with discretization
techniques not suitable for Total Lagrangian Formulation or Updated Lagrangian
Formulation.  In this paper the CR formulation is applied with Hybrid Trefftz Stress
(HTS) finite elements. Application of the CR formulation together with
HTS elements results in a powerful synergy, which is exemplified on a case study problem
presented in the paper. Although the presented methodology is applied to problems
of solid mechanics, the presentede method can be applied to image recovery,
where extraction of rigid body motion from image is an issue.

This paper follows earlier contributions \cite{Rankin,Crisfield1,Crisfield2,Felippa}.
Full literature review can be fond in \cite{Felippa}. Despite the large body of
literature on co-rotational formulations, formulation of a best fit rotation functional
remains an open reasearch question. In this paper we introduce a new general best-fit
rotator functional. Additionally a new Spin-Filter operator, mapping increment of degrees
of freedom into increment of rotation vector, consistent with new best-fit rotator
functional is determined.

\begin{figure}[!htp]
\centering
\includegraphics[width=0.9\textwidth]{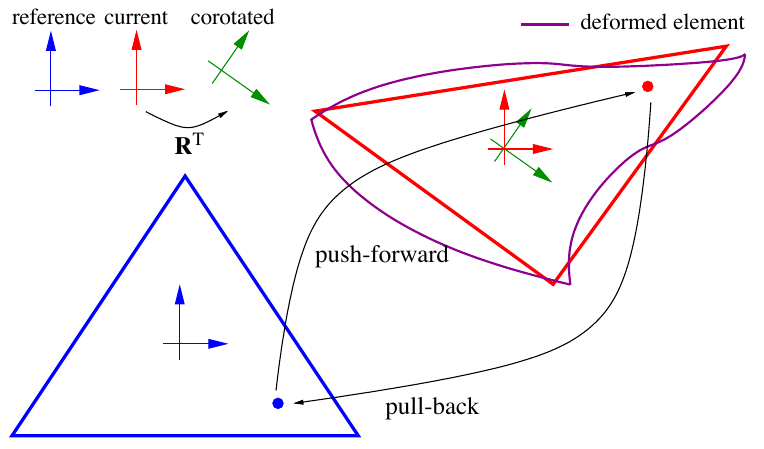}
\caption{\label{sec:HTS} Refrence and coorotated frames}
\end{figure}

In the first section approximation functions for the HTS elements are briefly
described. Next, kinematics of rotating deformable motion is described and the
best-firt rotator functional is formulated. Section \ref{sec:bestrorarir} is supplemented with
linearized equations, used to determine the rotation operator. In Section \ref{sec:proj}
a Spin-Liver linear operator transforming an increment of rotation angle into an increment
of displacement degrees of freedom, and Spin-Filter are formulated. With projection matrices at hand,
 a linearization equations of HTS elements are shown in Section \ref{sec:HTS}. In Section \ref{sec:examples} two
numerical examples are presented. The paper is concluded in Section \ref{sec:concl}.

\section{Approximation functions}

Lagrange CR formulation can be applied to a large class of discretization methods. In
this paper we apply it to a particular instance of the Hybrid Trefftz Stress elements.
In authors view combination of the HTS element with the CR formulation is a good example of
the robustness of the co-rotational approach: we apply the CR formulation to
a discontinuous displacement approximation, exclusively determined only on
element boundary. Moreover, unlike for the classical displacement Finite Element Method,
the employed displacement approximation does represent a partition of unity.

\subsection{Approximation of displacements on element faces}

\begin{figure}[!htp]
\centering
\includegraphics{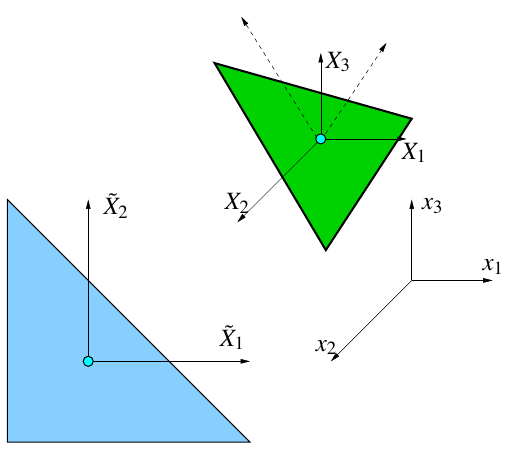}

\caption{Coordinate system on the reference element face and the element face in the reference configuration in space}

\end{figure}

In the spirit of hybrid FE formulations, displacement field is approximated on
element boundary, independently on each element face. In the presented approach we used
simple monomials given in the base coordinate system of element faces, see
Eq.~(\ref{eq:ugamma}). Displacements of face $f$ at point $\mathbf{X}_f$ (expressed in face coordinates) read
\begin{equation}
\mathbf{u}(\tilde{\mathbf{X}}_f) = 
\mathbf{U}_\Gamma(\tilde{\mathbf{X}}_f)\mathbf{q}. 
\end{equation}
Matrix of displacement approximation functions has simple form
\begin{equation} \label{eq:ugamma} 
\mathbf{U}_\Gamma(\tilde{\mathbf{X}}_f) = \left[
\begin{array}{ccc|cccccc|ccccccc} 1 & 0 & 0 & \tilde{X}_1 & 0 & 0 & \tilde{X}_2
& 0 & 0 & \tilde{X}^2_1 & 0 & 0 & \dots\\ 0 & 1 & 0 & 0 & \tilde{X}_1 & 0 & 0 &
\tilde{X}_2 & 0 & 0 & \tilde{X}^2_1 & 0 & \dots\\ 0 & 0 & 1 & 0 & 0 &
\tilde{X}_1 & 0 & 0 & \tilde{X}_2 & 0 & 0 & \tilde{X}^2_1 & \dots \end{array}
\right] 
\end{equation} where $\tilde{\mathbf{X}}$ is a coordinate vector in the reference system of element face.

Such approximation technique allows for an easy extension of the approximation space by arbitrary
functions, which are not necessarily a partition of unity. Moreover, since degrees of freedom are associated
with elements faces, for three-dimensional tetrahedrons (utilized here), an element face can have no more
than six neighbours (adjacent through elements on both sides). As a result, this formulation is suitable
for parallel computing where intercommunication between processors affects efficiency.

\subsection{Approximation of stresses in element volume}

Detailed description of HTS elements can be found in \cite{fraitas,Ktefftz}. In
this paper we limit ourselves to basic equations. In this section approximation
of stresses is briefly described. The Cauchy stress field within an element is approximated directly as:

\begin{equation}\label{eq:f1}
{\pmb \sigma} = \mathbf{S_v}(\mathbf{X}_f) \mathbf{v}
\end{equation}
where $\mathbf{S_v}$ is a matrix of field approximation functions and
$\mathbf{v}$ is the unknown vector of generalised stress degrees of freedom. We
note that $\mathbf{S_v}(\mathbf{X}_f)$ is defined on element face in Lagrange (co-rotated)
coordinates. In eq.~(\ref{eq:f1}), the stress approximation field is
chosen in order to automatically satisfy the equilibrium condition
(\ref{eq1}):
\begin{equation} \label{eq:f2}
\mathbf{L}^{T} \mathbf{S_v} = \mathbf{0}.
\end{equation}

where $\mathbf{L}$ is a differential operator related to the equilibrium equations.
According to the Cauchy equilibrium condition, the corresponding
stress induced traction field on the element boundary is
\begin{equation}
\label{eq:f7} \mathbf{t} =  \mathbf{n}{\pmb \sigma} =
\mathbf{nS_vv}
\end{equation}
Furthermore, assuming small strains and elastic body, the displacement field
within the domain is directly associated with this stress approximation
expressed as:
\begin{equation}
\label{eq:f7a} \mathbf{u}  =  \mathbf{U_vv}.
\end{equation}
where $\mathbf{U_v}$ is a displacement approximation matrix, such that
$\mathbf{S_v} = \mathbf{C}^{-1}\mathbf{L}\mathbf{U_v}$, where
$\mathbf{C}$ is a compliance matrix.

The stress field and the associated displacement field are derived from
appropriate Trefftz functions. For example, the stress and displacement field
could be given as polynomial functions derived from the Airy's stress function.
Alternatively, the Trefftz functions can be extended to functions which can
take into account singularities or are suitable for a particular boundary value
problem. Moreover, the approximation of stress is not {\em a priori} dependent
on a particular constitutive law. 

\section{Kinematics of CR frame}

For completeness some basic equations and theory of rotations are presented. 

Consider a body under rotational motion
\begin{equation}
\mathbf{x}(\mathbf{X}) = \mathbf{R}\mathbf{X}
\end{equation}
where $\mathbf{R}$ is $3$x$3$ is the rotation operator (also called rotor), $\mathbf{X}$ is a selected
material point and $\mathbf{x}$ is the spatial image of $\mathbf{X}$. Since $\mathbf{R}$ is the rotation
operator, distances between points and orientation are preserved, i.e.
$\mathbf{R}\mathbf{R}^\textrm{T} = \mathbf{I}$ and
$\textrm{det}(\mathbf{R})=1$.

In mathematical terms the set of all $3\times3$ orthogonal matrices with positive determinant forms the special
orthogonal group $SO(3)$, which is a manifold. Since rotations are points of the manifold $SO(3)$, at each
point $\mathbf{R}$ in $SO(3)$ a tangent Euclidean space $T_\mathbf{R} SO(3)$ collects rotation velocities $\dot{\mathbf{R}}$.

Consider now a nine-dimensional space of all $3\times3$ matrices and a surface defined by the constraint 
$f(\mathbf{R})=\mathbf{R}\mathbf{R}^\textrm{T}-\mathbf{I}=\mathbf{0}$.  The surface normal is orthogonal to vectors
 $\dot{\mathbf{R}} \in T_\mathbf{R} SO(3)$ or in other words

\begin{equation} \label{eq:k2}
\dot{\mathbf{R}}\mathbf{R}^\textrm{T}  + \mathbf{R}\dot{\mathbf{R}}^\textrm{T} = 0.
\end{equation}

Let us define ${\boldsymbol \omega}$ as
\begin{equation} \label{eq:spin}
{\boldsymbol \omega} = \dot{\mathbf{R}}\mathbf{R}^\textrm{T}.
\end{equation}
Since equation (\ref{eq:k2}) is satisfied for all $\mathbf{R} \in SO(3)$ and
$\dot{\mathbf{R}} \in T_\mathbf{R} SO(3)$ we note that  $\boldsymbol \omega$ is anti-symmetric, ${\boldsymbol \omega}^\textrm{T}=-{\boldsymbol \omega}$.

The operator $\boldsymbol \omega$ bears the name of the spatial angular velocity
\begin{equation} \label{eq:velocity_ref}
\dot{\mathbf{x}} = {\boldsymbol \omega}\mathbf{x}.
\end{equation}

For $\boldsymbol \omega = \textrm{const}$ the ordinary differential equation (\ref{eq:velocity_ref}) has a known solution
\begin{equation}
\mathbf{x}(t) = \textrm{Exp}(t{\boldsymbol \omega})\mathbf{R}(0)\mathbf{X}(0).
\end{equation}
where time is expressed explicitly by $t$ and $\textrm{Exp}$ stands for the matrix exponential.
Although in general $\boldsymbol \omega$ is not constant in time, the above formula serves as a useful
device for implementing incremental updates of $\mathbf{R}$.

\subsection{Variations of rotation operator}

With the brief theory at hand, we introduce variations in order to construct
tangent operators for the HTS formulation.  From (\ref{eq:spin}) we can see that
$\dot{\mathbf{R}} = {\boldsymbol \omega}\mathbf{R}$ and hence the variation of $\mathbf{R}$
can be expressed in the spatial form

\begin{equation} \label{eq:var_rot}
\delta \mathbf{R} = \delta {\boldsymbol \Phi} \mathbf{R}
\end{equation}
where $\delta {\boldsymbol \Phi} \in T_\mathbf{R} SO(3)$.

\subsection{Variation of displacements}

For simplicity and without loss of generality we analyse motion
without the rigid body translations. New position of a particle is given by
\begin{equation}
\mathbf{x}^r(\mathbf{X}) = \mathbf{R}\mathbf{X}
\end{equation}
where $\mathbf{X}$ is the position vector in the reference configuration. With that at
hand we define displacement by
\begin{equation}
\mathbf{u}^r(\mathbf{X}) = \mathbf{x}^r(\mathbf{X})-\mathbf{X} = (\mathbf{R}-\mathbf{I})\mathbf{X}
\end{equation}

Displacements are additively decomposed into the deformational and the rotational part, which yields
\begin{equation} \label{eq:ud}
\mathbf{u}^d(\mathbf{X}) = \mathbf{u}(\mathbf{X}) + (\mathbf{I}-\mathbf{R})\mathbf{X} =
\mathbf{u}(\mathbf{X}) - \mathbf{u}^r(\mathbf{X}) 
\end{equation}
where $\mathbf{u}^\textrm{d}$ and $\mathbf{u}^r$ are the deformational and
the rotational parts of displacements, respectively.

In order to derive tangent operators, a useful variation of deformational displacements (\ref{eq:ud}) is given by
\begin{equation}
\begin{split}
\delta \mathbf{u}^d(\mathbf{X}) &=
\delta \mathbf{u} - \delta \mathbf{R} \mathbf{X} \\
&= \delta \mathbf{u} - \textrm{Spin}[\delta\vec{\boldsymbol\Phi}] \mathbf{R} (\mathbf{X} - \mathbf{a}) \\ 
&= \delta \mathbf{u} - \textrm{Spin}[\delta\vec{\boldsymbol\Phi}] \mathbf{x}^\textrm{r} \\
&= \delta \mathbf{u} + \textrm{Spin}[\mathbf{x}^r] \delta\vec{\boldsymbol\Phi}
\end{split}
\end{equation}
where variation of the axial rotation vector $\delta\vec{\boldsymbol\Phi}$ is
used. We note, that $\delta\vec{\boldsymbol\Phi}$ is a pseudo-vector, i.e. a 3-dimensional representation
of a $3\times3$ antisymmetric matrix. The operator $\textrm{Spin}[\cdot]$ takes form
\begin{equation}
\textrm{Spin}[\mathbf{x}] = 
\left[
\begin{array}{ccc}
0 & -x_3 & x_2 \\
x_3 & 0 & -x_1 \\
-x^2 & x_1 & 0
\end{array}
\right]
\end{equation} 
for an arbitrary vector $\mathbf{x}$.

\section{Best fit rotor \label{sec:bestrorarir}} 

Computation of the rotation operator for a given displacement field is at the hart of
the co-rotational formulation. In our view, an ideal co-rotational formulation should address the following objectives:

\begin{itemize}
\item \emph{Versatility.} It must work for static and dynamic problems, with arbitrary  discretization of displacements.
\item \emph{Rigid bodies.} For rigid body motion, it should give the same results as a body-attached frame.
This simplifies coupling of FEM-CR methods with multibody dynamics codes \cite{Felippa}.
\item \emph{Finite stretch and pure sheer.} If a element is subjected to a uniform
finite stretch or pure sheer, no spurious rotation should appear.
\end{itemize}

In order to find the rotation operator $\mathbf{R}$, a new functional is proposed here.
Stationary points of this functional determine rotation, which complies with the above objectives.

We note, that for a rotation free deformation, the antisymmetric part of the displacement
gradient disappears. With that observation at hand, assuming that the rotation is
constant in the element domain, the antisymmetric part of the volume averaged deformational displacement
is used to construct the best fit rotator.
\begin{equation}
\int_\Omega 
\left(
\left[\frac{\partial\mathbf{u}^\textrm{d}}{\partial \mathbf{X}}\right] - 
\left[\frac{\partial\mathbf{u}^\textrm{d}}{\partial \mathbf{X}}\right]^\textrm{T}
\right)
\textrm{d}\Omega = \mathbf{0}
\end{equation}

Noting that the antisymmetric part of the displacement gradient has three non-zero
linearly independent components, the antisymmetric part of the displacement gradient can
be uniquely expressed  by a pseudo-vector $\vec{\mathbf{h}}$. Utilizing that and
applying integration by parts, a pseudo-vector $\vec{\mathbf{h}}$ can be expressed by 
the boundary integral
\begin{equation}
\vec{\mathbf{h}} = 
\int_\Gamma
\mathbf{n} \times \mathbf{u}^d
\textrm{d}\Gamma
=
\int_\Gamma
\textrm{Spin}[\mathbf{n}]\mathbf{u}^d
\textrm{d}\Gamma
\end{equation}
where $\mathbf{n}$ is the spatial normal vector field on the element boundary $\Gamma$ in the current
or reference configuration (since the motion is roughly rigid we can assume the Jacobian to be nearly one).
A length of the pseudo-vector $\vec{\mathbf{h}}$ is used to formulate the best fit functional. Utilizing that $\int_\Gamma
\textrm{Spin}[\mathbf{N}]\mathbf{X}\textrm{d}\Gamma = 0$, the functional takes form
\begin{equation} \label{eq:bestfit}
\begin{split}
\mathbf{J}(\mathbf{R}) &= 
\frac{1}{2}
[\vec{\mathbf{h}}]^\textrm{T}[\vec{\mathbf{h}}]
= \frac{1}{2}
\left[
\int_\Gamma
\textrm{Spin}[\mathbf{n}]\mathbf{u}^d
\textrm{d}\Gamma
\right]^\textrm{T}
\left[
\int_\Gamma
\textrm{Spin}[\mathbf{n}]\mathbf{u}^d
\textrm{d}\Gamma
\right] \\ &= \frac{1}{2}
\left[
\int_\Gamma
\textrm{Spin}[\mathbf{R}\mathbf{N}](\mathbf{u}-\mathbf{u}^\textrm{r})
\textrm{d}\Gamma
\right]^\textrm{T}
\left[
\int_\Gamma
\textrm{Spin}[\mathbf{R}\mathbf{N}](\mathbf{u}-\mathbf{u}^\textrm{r})
\textrm{d}\Gamma
\right] \\ &= \frac{1}{2}
\left[
\int_\Gamma
\textrm{Spin}[\mathbf{N}]\mathbf{R}^\textrm{T}\mathbf{x}
\textrm{d}\Gamma
\right]^\textrm{T}
\left[
\int_\Gamma
\textrm{Spin}[\mathbf{N}]\mathbf{R}^\textrm{T}\mathbf{x}
\textrm{d}\Gamma
\right]
\end{split}
\end{equation}
where $\mathbf{N}$ is the referential normal vector field on the element boundary $\Gamma$ in the current
or reference configuration. We note, that the functional expressed in the above form can be applied
to problems where displacements are only known on elements boundaries (and can be discontinuous).

For given displacements $\mathbf{u}$, a new position vector $\mathbf{x}$ in the current
configuration is easy to determine. For a fixed $\mathbf{x}$, the condition of stationarity 
of $J(\mathbf{R})$ represents the set of nonlinear equations for finding $\mathbf{R}$
\begin{equation}
\begin{split}
dJ &= \frac{\partial J(\mathbf{R})}{\partial\mathbf{R}} : \delta\mathbf{R} \\ &=
\left[ \int_\Gamma
\textrm{Spin}[\mathbf{N}]\mathbf{R}^\textrm{T}\mathbf{x}
\textrm{d}\Gamma \right]^\textrm{T}
\left[ \int_\Gamma
\textrm{Spin}[\mathbf{N}]\delta\mathbf{R}^\textrm{T}\mathbf{x}
\textrm{d}\Gamma \right] \\ &=
\left[ \int_\Gamma
\textrm{Spin}[\mathbf{N}]\mathbf{R}^\textrm{T}\mathbf{x}
\textrm{d}\Gamma \right]^\textrm{T}
\left[ \int_\Gamma
\textrm{Spin}[\mathbf{N}]
\mathbf{R}^\textrm{T}\textrm{Spin}[\mathbf{x}]
\textrm{d}\Gamma
\right]\delta\vec{\boldsymbol\Phi} = 0.
\end{split}
\end{equation}
where (\ref{eq:var_rot}) has been utilized. We then requre that
\begin{equation}
\begin{split}
\mathbf{H} &=
\left[ \int_\Gamma
\textrm{Spin}[\mathbf{N}]\mathbf{R}^\textrm{T}\mathbf{x}
\textrm{d}\Gamma \right]^\textrm{T}
\left[ \int_\Gamma
\textrm{Spin}[\mathbf{N}]
\mathbf{R}^\textrm{T}\textrm{Spin}[\mathbf{x}]
\textrm{d}\Gamma
\right] = \mathbf{0}
\end{split}
\end{equation}
and solve the above sytem by means of the Newton-Raphson iterations
\begin{equation} \label{eq:newton0}
\frac{\partial \mathbf{H}}{\partial {\boldsymbol \Phi}}\Delta {\boldsymbol \Phi} = -\mathbf{H}
\end{equation}
\begin{equation} \label{eq:newton1}
\mathbf{R}_\textrm{k+1} =
\textrm{Exp}(\Delta{\boldsymbol \Phi})
\mathbf{R}_\textrm{k}.
\end{equation}

In order to approximate (\ref{eq:newton0})
we linearise $\textrm{Exp}(\Delta {\boldsymbol \Phi})$
\begin{equation}
\begin{split}
\textrm{Exp}\left(\Delta {\boldsymbol \Phi}\right) &= 
\mathbf{I} + 
\frac
{\sin  \parallel \Delta \vec{\boldsymbol \Phi} \parallel}
{\parallel \Delta \vec{\boldsymbol \Phi} \parallel}
{\Delta {\boldsymbol \Phi}}+
\frac{1 - \cos  \parallel \vec{\Delta {\boldsymbol \Phi}} \parallel}
{ \parallel \Delta \vec{\boldsymbol \Phi} \parallel^2}
{\Delta {\boldsymbol \Phi}}^2\\
&\approx
\mathbf{I} + {\Delta {\boldsymbol \Phi}}
\end{split}
\end{equation}
and noting that
$\mathbf{R}_\textrm{k+1} \approx \mathbf{R}_\textrm{k}(\mathbf{I}+\Delta {\boldsymbol \Phi})$
we represent (\ref{eq:newton0}) by up to first order terms of
$\mathbf{H}(\mathbf{R}_\textrm{k}(\mathbf{I}+\Delta {\boldsymbol \Phi})) = \mathbf{0}$. That is

\begin{equation}
\begin{split}
&\left[
\int_\Gamma
\textrm{Spin}[\mathbf{N}]
\mathbf{R}_k^\textrm{T}\textrm{Spin}[\mathbf{x}]
\textrm{d}\Gamma
\right]^\textrm{T}
\left[
\int_\Gamma
\textrm{Spin}[\mathbf{N}]\mathbf{R}_k^\textrm{T}\textrm{Spin}[\mathbf{x}]
\textrm{d}\Gamma
\right]\Delta\vec{\boldsymbol \Phi}
\\ &+
\mathbf{e}_i\otimes
\left[
\int_\Gamma
\textrm{Spin}[\mathbf{N}]\mathbf{R}_k^\textrm{T}\mathbf{x}
\textrm{d}\Gamma
\right]^\textrm{T}
\left[
\int_\Gamma
\textrm{Spin}[\mathbf{N}]
\mathbf{R}_k^\textrm{T}\textrm{Spin}\left[\textrm{Spin}_j[\mathbf{x}]\right]
\textrm{d}\Gamma
\right]
\Delta\vec{\boldsymbol \Phi} \\ &=
-
\left[
\int_\Gamma
\textrm{Spin}[\mathbf{N}]\mathbf{R}_k^\textrm{T}\mathbf{x}
\textrm{d}\Gamma
\right]^\textrm{T} 
\left[
\int_\Gamma
\textrm{Spin}[\mathbf{N}]
\mathbf{R}_k^\textrm{T}\textrm{Spin}[\mathbf{x}]
\textrm{d}\Gamma
\right]
\end{split}
\end{equation}
where $\mathbf{e}_i$ is the Cartesian base vector, $i,j \in
\left\{1,2,3\right\}$ and $\textrm{Spin}_j[\mathbf{x}]$ represents the $j$th
column of $\textrm{Spin}[\mathbf{x}]$.

In this section a new best-fit functional, based on boundary integrals is
presented. Stationary point of this functional lead to a system of non-linear
algebraic equations, form which current rotation operator can be calculated by means
of Newton iterations. It should be noted, that the above procedure allows to determine
the rotation quite independently of the motion history. However, in order to avoid troubles
with multiple minima and correctly trace the history of rotation it is appropriate to
initialize the iterative procedure with the most recent value of rotation.

\section{Projection matrices} \label{sec:proj}

In order to construct tangent operators for the consistent Newton scheme projection
matrices filtering a variation of the deformational displacements from the variation
of the total displacements have to be derived.
\begin{equation}
\delta \mathbf{q}^\textrm{d} = \mathbf{P}\delta\mathbf{q} = 
(\mathbf{I}-\mathbf{SG})\delta\mathbf{q}
\end{equation}
Projection matrix $\mathbf{P}$, Spin-Liver $\mathbf{S}$ and Spin-Filter
$\mathbf{G}$ depend on the approximation method. In this paper we construct
these matrices for the particular approximation of HTS finite elements. However,
the methodology is general and can be applied to a large class approximation methods.


\subsubsection{Spin Liver}

\begin{figure}[!htp]
\centering
\includegraphics[width=0.99\textwidth]{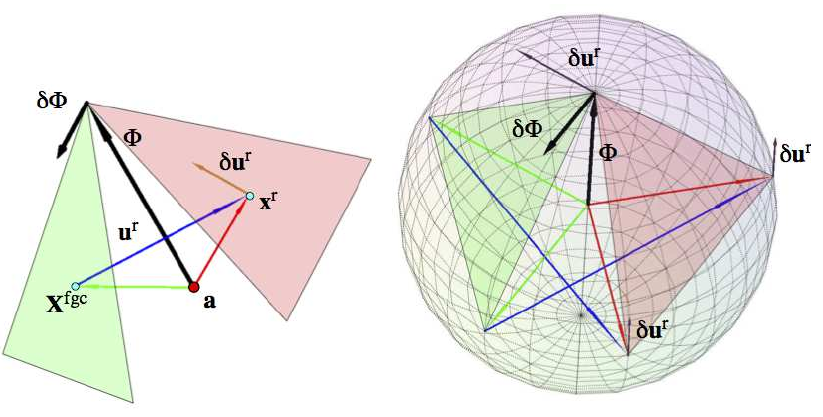}

\caption{For a given current face frame and an arbitrary variation of spin
$\delta\boldsymbol\Omega$, the variation of rotational linear displacements are
presented in figure on the left-hand side. Linear rotational displacements are
computed using formula $\delta\mathbf{u}^\textrm{r} =
\mathbf{S}\delta\vec{\boldsymbol\Phi}$. The green triangular face represents
the reference configuration, the red triangular face depicts the current configuration.
The Spin-Liver is constructed so that the resulting variations of the linear displacements are
orthogonal to the variation of the rotation axis and tangent to the sphere, see figure
on the right-hand side. Position of points on the sphere depends on the current
rotation axis and rotation angle.}

\end{figure}

A Spin-Liver matrix represents a linear operator transforming a variation of
the rotation axis into a variation of linear displacements. Spin-Liver $\mathbf{S}$
depends only on constant and linear approximation functions. Higher order
terms in displacements approximation basis are neglected in the formulation of
Spin-Liver. It can be noted, that for any displacement function given by
higher-order terms the motion is not stress free, hence higher-order therms cannot
be a part of the rigid body motion. 

Exploiting features of the HTS displacement approximation, see Eq.~(\ref{eq:ugamma}), the
Spin-Liver $\mathbf{S}$ can be defined for each face independently, which yields 
\begin{equation} 
\delta\mathbf{u}^r(\mathbf{X}) 
= \mathbf{U}_\Gamma(\tilde{\mathbf{X}}_f)
\mathbf{S}_f(\mathbf{X}^\textrm{fgc})\delta \vec{\boldsymbol \Omega} 
\end{equation} 
where $\mathbf{X}^\textrm{fgc}=\mathbf{X}(\tilde{\mathbf{X}}=\mathbf{0})$ is
a face geometrical center. The Spin-Liver matrix itself, for our particular implementation
of the HTS element has simple form
\begin{equation}
\mathbf{S}(\mathbf{X}^\textrm{fgc})= 
\left[
\begin{array}{ccc}
0 & x_3^r & -x_2^r \\
-x_3^r & 0 & x_1^r \\
x_2^r & -x_1^r & 0 \\
\hline
0 & \partial x_3^r/\partial \tilde{X}_1 & -\partial x_2^r/\partial \tilde{X}_1 \\
-\partial x_3^r/\partial \tilde{X}_1 & 0 & \partial x_1^r/\partial \tilde{X}_1 \\
\partial x_2^r/\partial \tilde{X}_1 & -\partial x_1^r/\partial \tilde{X}_1 & 0 \\
0 & \partial x_3^r/\partial \tilde{X}_2 & -\partial x_2^r/\partial \tilde{X}_2 \\
-\partial x_3^r/\partial \tilde{X}_2 & 0 & \partial x_1^r/\partial \tilde{X}_2 \\
\partial x_2^r/\partial \tilde{X}_2 & -\partial x_1^r/\partial \tilde{X}_2 & 0 \\
\hline
0 & 0 & 0 \\
\dots
\end{array}
\right]
\end{equation}
where $\mathbf{x}^\textrm{r} = \mathbf{R}\mathbf{X}^\textrm{fgc}$ and 
$
\partial \mathbf{x}^r/\partial \tilde{\mathbf{X}}|_{\mathbf{X}^\textrm{fgc}} = 
\mathbf{R} \partial \mathbf{X}/\partial \tilde{\mathbf{X}}|_{\mathbf{X}^\textrm{fgc}}
$. Definition of the Spin-Liver for other types of elements can be found in \cite{Felippa}.

\subsubsection{Spin Filter}

\begin{figure}[!htp]
\centering
\includegraphics[width=0.9\textwidth]{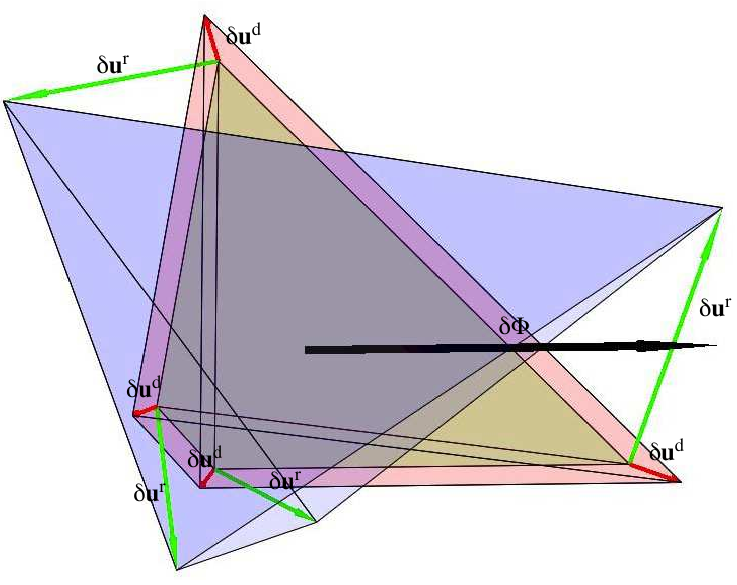}

\caption{ 
Current configuration is given by the green tetrahedron. Variation of current
configuration is depicted by the blue tetrahedron. Variation of the rotation axis is
visualized by the black arrow. Variation of the deformation is depicted by the red
tetrahedral.
}

\end{figure}

Matrix $\mathbf{G}$ represents a linear operator which transforms a
variation of linear displacements into a variation of rotation axis. In other words
Spin-Filter accounts for a change of magnitude and direction of the rotation
axis due to the change of linear displacements.

In order to get quadratic convergence, Spin-Filter matrix has to be consistent
with functional $J$, defined in Section \ref{sec:bestrorarir}.  For a
given current configuration, i.e. rotation $\mathbf{R}$, and linearized
displacements, functional $J$ takes form 
\begin{equation}
\begin{split}
&J(\delta\vec{\boldsymbol\Phi}) \\ &= \frac{1}{2}
\left[
\int_\Gamma
\textrm{Spin}[\mathbf{N}]\mathbf{R}^\textrm{T}(\delta\mathbf{u}+\textrm{Spin}[\mathbf{x}^\textrm{r}]\delta\vec{\boldsymbol\Phi})
\textrm{d}\Gamma
\right]^\textrm{T}
\left[
\int_\Gamma
\textrm{Spin}[\mathbf{N}]\mathbf{R}^\textrm{T}(\delta\mathbf{u}+\textrm{Spin}[\mathbf{x}^\textrm{r}]\delta\vec{\boldsymbol\Phi})
\textrm{d}\Gamma
\right]
\end{split}
\end{equation}
Stationary point of this functional yields 
\begin{equation}
\begin{split}
&\frac{\partial J(\delta\vec{\boldsymbol\Phi})}{\partial\delta\vec{\boldsymbol\Phi}} : \delta\vec{\boldsymbol\Phi} = \\
&\left[
\int_\Gamma
\textrm{Spin}[\mathbf{N}]\mathbf{R}^\textrm{T}\textrm{Spin}[\mathbf{x}^\textrm{r}]
\textrm{d}\Gamma
\right]^\textrm{T}
\left[
\int_\Gamma
\textrm{Spin}[\mathbf{N}]\mathbf{R}^\textrm{T}(\delta\mathbf{u}+\textrm{Spin}[\mathbf{x}^\textrm{r}]\delta\vec{\boldsymbol\Phi})
\textrm{d}\Gamma
\right]\\ &= 0.
\end{split}
\end{equation}
Equation of the stationary point is given by the system of three linear equations 
\begin{equation}
\begin{split}
\left[
\int_\Gamma
\textrm{Spin}[\mathbf{N}]\mathbf{R}^\textrm{T}\textrm{Spin}[\mathbf{x}^\textrm{r}]
\textrm{d}\Gamma
\right]^\textrm{T}
\left[
\int_\Gamma
\textrm{Spin}[\mathbf{N}]\mathbf{R}^\textrm{T}\textrm{Spin}[\mathbf{x}^\textrm{r}]
\textrm{d}\Gamma
\right] 
\delta\vec{\boldsymbol\Phi}
= \\
-\left[
\int_\Gamma
\textrm{Spin}[\mathbf{N}]\mathbf{R}^\textrm{T}\textrm{Spin}[\mathbf{x}^\textrm{r}]
\textrm{d}\Gamma
\right]^\textrm{T}
\left[
\int_\Gamma
\textrm{Spin}[\mathbf{N}]\mathbf{R}^\textrm{T}\delta\mathbf{u}
\textrm{d}\Gamma
\right]. 
\end{split}
\end{equation}
Substituting $\delta\mathbf{u} = \mathbf{U}_\Gamma \delta\mathbf{q}$ yields
\begin{equation}
\begin{split}
\left[
\int_\Gamma
\textrm{Spin}[\mathbf{N}]\mathbf{R}^\textrm{T}\textrm{Spin}[\mathbf{x}^\textrm{r}]
\textrm{d}\Gamma
\right]^\textrm{T}
\left[
\int_\Gamma
\textrm{Spin}[\mathbf{N}]\mathbf{R}^\textrm{T}\textrm{Spin}[\mathbf{x}^\textrm{r}]
\textrm{d}\Gamma
\right] 
\delta\vec{\boldsymbol\Phi}
= \\
-\left[
\int_\Gamma
\textrm{Spin}[\mathbf{N}]\mathbf{R}^\textrm{T}\textrm{Spin}[\mathbf{x}^\textrm{r}]
\textrm{d}\Gamma
\right]^\textrm{T}
\left[
\int_\Gamma
\textrm{Spin}[\mathbf{N}]\mathbf{R}^\textrm{T}\mathbf{U}_\Gamma
\textrm{d}\Gamma
\right] 
\delta\mathbf{q}
\end{split}
\end{equation}

As a result, a matrix representing the Spin-Filter
\begin{equation} \label{eq:sp1}
\delta\vec{\boldsymbol\Phi} = \mathbf{G}\delta\mathbf{q},
\end{equation}
is expressed in the form
\begin{equation} \label{eq:sp2}
\begin{split}
\mathbf{G} = 
-\left\{
\left[
\int_\Gamma
\textrm{Spin}[\mathbf{N}]\mathbf{R}^\textrm{T}\textrm{Spin}[\mathbf{x}^\textrm{r}]
\textrm{d}\Gamma
\right]^\textrm{T}
\left[
\int_\Gamma
\textrm{Spin}[\mathbf{N}]\mathbf{R}^\textrm{T}\textrm{Spin}[\mathbf{x}^\textrm{r}]
\textrm{d}\Gamma
\right]
\right\}^{-1}\\
\left[
\int_\Gamma
\textrm{Spin}[\mathbf{N}]\mathbf{R}^\textrm{T}\textrm{Spin}[\mathbf{x}^\textrm{r}]
\textrm{d}\Gamma
\right]^\textrm{T}
\left[
\int_\Gamma
\textrm{Spin}[\mathbf{N}]\mathbf{R}^\textrm{T}\mathbf{U}_\Gamma
\textrm{d}\Gamma
\right].
\end{split}
\end{equation}

The novelty of this approach in comparison to those found in other papers related to the CR
formulation, is that the Spin-Filter is derived form best fit functional and can be
therefore applied to a large class of approximation methods.

\section{HTS formulation} 

Although the co-rotational formulation can be derived independently form a problem
formulation by consistent differentiation of projection matrices
\cite{Felippa}, in this paper we will formulate tangent operators and residuals
directly from weak forms. In authors opining this allows for better
understanding of the problem, without loss of generality. 

The co-rotational Hybrid Trefftz Stress element is built using 
Lagrange coordinates, i.e. integrals are over reference volumes and faces.
Stresses are expressed in the co-rotated (current) configuration and
displacements are expressed in the reference coordinate system. At each iteration
co-rotated (Cauchy) stress is pulled-back to the reference coordinate system (or
conversly, the displacement is pushed-forward to the co-rotated reference system).

\subsection{Integral form}

Static boundary conditions are enforced in strong integral form, i.e. we weight
static boundary condition by function $\mathbf{w_1}$ and integrate over element
boundary. 
\begin{equation}\label{eq:hts1}  
{\mathbf{R}_\sigma} 
=\int_{\Gamma}\mathbf{w_1}^\textrm{T}\mathbf{R}\boldsymbol{\hat{\sigma}}\mathbf{\hat{n}}\textrm{d}\Gamma 
-\int_{\Gamma}\mathbf{w_1}^\textrm{T}\overline{\mathbf{t}}\textrm{d}\Gamma
= \mathbf{0}
\end{equation}
Function $\mathbf{w_1}$ is given in the reference coordinate system and stresses
are pulled-back to the reference coordinate system. This results in geometrical
nonlinearities, i.e. the integral becomes a function of the unknown rotor $\mathbf{R}$.

\subsection{Weak form}

The compatibly conditions within a body domain are imposed in a weak sense, where by
the differentiation by parts Dirichlet boundary conditions are included 
\begin{equation}\label{eq:Res_epsilon} 
\begin{split}
{\mathbf{R}_\varepsilon} 
=\int_{\Omega}\textrm{tr}[\mathbf{w_2}\boldsymbol{\varepsilon}]\textrm{d}\Omega 
-\int_{\Gamma_\textrm{u}}(\mathbf{w_2n})^\textrm{T}(\mathbf{u}_\Gamma-\mathbf{c}-\mathbf{u}^r)\textrm{d}\Gamma
-\int_{\Gamma_\sigma}(\mathbf{w_2n})^\textrm{T}\overline{\mathbf{u}}\textrm{d}\Gamma 
= \mathbf{0}
\end{split}
\end{equation}
\begin{equation}
\begin{split}
{\mathbf{R}_\varepsilon} 
=\int_{\Omega}\textrm{tr}[\mathbf{w_2}\boldsymbol{\varepsilon}]\textrm{d}\Omega 
-\int_{\Gamma_\textrm{u}}(\mathbf{w_2n})^\textrm{T}\mathbf{u}^\textrm{d}_\Gamma\textrm{d}\Gamma
-\int_{\Gamma_\sigma}(\mathbf{w_2n})^\textrm{T}\overline{\mathbf{u}}\textrm{d}\Gamma 
= \mathbf{0}
\end{split}
\end{equation}
\begin{equation}
\begin{split}
{\mathbf{R}_\varepsilon} 
=\int_{\Omega}\textrm{tr}[\mathbf{R}\mathbf{\hat{w}_2}\mathbf{R}^\textrm{T}\boldsymbol{\varepsilon}]\textrm{d}\Omega 
-\int_{\Gamma_\textrm{u}}(\mathbf{R}\mathbf{\hat{w}_2\hat{n}})^\textrm{T}\mathbf{u}^\textrm{d}_\Gamma\textrm{d}\Gamma
-\int_{\Gamma_\sigma}(\mathbf{R}\mathbf{\hat{w}_2\hat{n}})^\textrm{T}\overline{\mathbf{u}}\textrm{d}\Gamma 
= \mathbf{0}
\end{split}
\end{equation}
\begin{equation}
\begin{split}
{\mathbf{R}_\varepsilon} 
=\int_{\Omega}\textrm{tr}[\mathbf{R}\mathbf{\hat{w}_2}\mathbf{R}^\textrm{T}\mathbf{R}\boldsymbol{\hat{\varepsilon}}\mathbf{R}^\textrm{T}]\textrm{d}\Omega 
-\int_{\Gamma_\textrm{u}}(\mathbf{R}\mathbf{\hat{w}_2\hat{n}})^\textrm{T}\mathbf{u}^\textrm{d}_\Gamma\textrm{d}\Gamma
-\int_{\Gamma_\sigma}(\mathbf{R}\mathbf{\hat{w}_2\hat{n}})^\textrm{T}\overline{\mathbf{u}}\textrm{d}\Gamma 
= \mathbf{0}
\end{split}
\end{equation}
\begin{equation} \label{eq:hts2}
\begin{split}
{\mathbf{R}_\varepsilon} 
=\int_{\Omega}\textrm{tr}[\mathbf{\hat{w}_2}\boldsymbol{\hat{\varepsilon}}]\textrm{d}\Omega 
-\int_{\Gamma_\textrm{u}}(\mathbf{\hat{w}_2\hat{n}})^\textrm{T}\mathbf{R}^\textrm{T}\mathbf{u}^\textrm{d}_\Gamma\textrm{d}\Gamma
-\int_{\Gamma_\sigma}(\mathbf{\hat{w}_2\hat{n}})^\textrm{T}\mathbf{R}^\textrm{T}\overline{\mathbf{u}}\textrm{d}\Gamma 
= \mathbf{0}
\end{split}
\end{equation}
A weighting function $\mathbf{w_2}$ is defined in the co-rotated coordinate system and is transformed
to the reference coordinate system. We note, that the term under the volume integral is a scalar, so it
not depend on the coordinate system. As a result it is constant for all steps analysis and is computed
only once.

\subsection{Linearization}

The system of nonlinear algebraic equations is linearized and solved using Newton
method. The rotation operator $\mathbf{R}$ is a function of current displacements
and it is computed at each iteration using the best-fit functional (\ref {eq:bestfit}). 
At a current cofiguratiom, linearized Eq.~(\ref{eq:hts1}) takes from
\begin{equation} 
{\mathbf{R}_\sigma} 
+\int_{\Gamma}\mathbf{w_1}^\textrm{T}\mathbf{R}\delta\boldsymbol{\hat{\sigma}}\mathbf{\hat{n}}\textrm{d}\Gamma 
+\int_{\Gamma}\mathbf{w_1}^\textrm{T}\delta\mathbf{R}\boldsymbol{\hat{\sigma}}\mathbf{\hat{n}}\textrm{d}\Gamma 
= \mathbf{0}
\end{equation}
\begin{equation} 
{\mathbf{R}_\sigma} 
+\int_{\Gamma}\mathbf{w_1}^\textrm{T}\mathbf{R}\delta\boldsymbol{\hat{\sigma}}\mathbf{\hat{n}}\textrm{d}\Gamma 
+\int_{\Gamma}\mathbf{w_1}^\textrm{T}\delta{\boldsymbol\Phi}\mathbf{R}\boldsymbol{\hat{\sigma}}\mathbf{\hat{n}}\textrm{d}\Gamma 
= \mathbf{0}
\end{equation}
\begin{equation} 
{\mathbf{R}_\sigma} 
+\int_{\Gamma}\mathbf{w_1}^\textrm{T}\mathbf{R}\delta\boldsymbol{\hat{\sigma}}\mathbf{\hat{n}}\textrm{d}\Gamma 
-\int_{\Gamma}\mathbf{w_1}^\textrm{T}\textrm{Spin}[\mathbf{R}\boldsymbol{\hat{\sigma}}\mathbf{\hat{n}}]\textrm{d}\Gamma 
\delta\vec{\boldsymbol\Phi}
= \mathbf{0}
\end{equation}

In a similar manner to Eq.~(\ref{eq:hts1}), residual Eq.~(\ref{eq:hts2}) is linearized and yields
\begin{equation}
\begin{split}
{\mathbf{R}_\varepsilon}
+\int_{\Omega}\textrm{tr}[\mathbf{\hat{w}_2}\delta\boldsymbol{\hat{\varepsilon}}]\textrm{d}\Omega 
-\int_{\Gamma}(\mathbf{\hat{w}_2\hat{n}})^\textrm{T}\mathbf{R}^\textrm{T}\delta\mathbf{u}^\textrm{d}_\Gamma\textrm{d}\Gamma
-\int_{\Gamma}(\mathbf{\hat{w}_2\hat{n}})^\textrm{T}\delta\mathbf{R}^\textrm{T}\mathbf{u}^\textrm{d}_\Gamma\textrm{d}\Gamma
= \mathbf{0}
\end{split}
\end{equation}
\begin{equation}
\begin{split}
{\mathbf{R}_\varepsilon}
+\int_{\Omega}\textrm{tr}[\mathbf{\hat{w}_2}\delta\boldsymbol{\hat{\varepsilon}}]\textrm{d}\Omega 
-\int_{\Gamma}(\mathbf{\hat{w}_2\hat{n}})^\textrm{T}\mathbf{R}^\textrm{T}\delta\mathbf{u}^\textrm{d}_\Gamma\textrm{d}\Gamma
+\int_{\Gamma}(\mathbf{\hat{w}_2\hat{n}})^\textrm{T}\mathbf{R}^\textrm{T}\delta{\boldsymbol\Phi}\mathbf{u}^\textrm{d}_\Gamma\textrm{d}\Gamma
= \mathbf{0}
\end{split}
\end{equation}
\begin{equation}
\begin{split}
{\mathbf{R}_\varepsilon}
+\int_{\Omega}\textrm{tr}[\mathbf{\hat{w}_2}\delta\boldsymbol{\hat{\varepsilon}}]\textrm{d}\Omega 
-\int_{\Gamma}(\mathbf{\hat{w}_2\hat{n}})^\textrm{T}\mathbf{R}^\textrm{T}\delta\mathbf{u}^\textrm{d}_\Gamma\textrm{d}\Gamma
-\int_{\Gamma}(\mathbf{\hat{w}_2\hat{n}})^\textrm{T}\mathbf{R}^\textrm{T}\textrm{Spin}[\mathbf{u}^\textrm{d}_\Gamma]\textrm{d}\Gamma
\delta\vec{\boldsymbol\Phi}
= \mathbf{0}
\end{split}
\end{equation}

A key term in above equations is the change of the angular vector
$\delta\vec{\boldsymbol\Phi}$. We note $\delta\vec{\boldsymbol\Phi}$ depends on
the current configuration and is a linear function of the variation of displacement
degrees of freedom, see~Eq.~(\ref{eq:sp1}) and Eq.~(\ref{eq:sp2}).

\subsection{Matrix form}

The linearized equations are reformulated, with stresses and displacements
substituted by suitable products approximation functions and increments of vectors
storing stress and displacement degrees of freedom. Spin-Filter $\mathbf{G}$ and Spin-Liver
$\mathbf{S}$ are utilized in order to express linearized equations exclusively
in terms unknown vectors, which yields
\begin{equation} \label{eq:res1}
{\mathbf{R}_\sigma}  
+\int_{\Gamma}\mathbf{U}_\Gamma^\textrm{T}\mathbf{\hat{n}\hat{S}_v}\textrm{d}\Gamma\delta\mathbf{v}
-\int_{\Gamma}\mathbf{U}_\Gamma^\textrm{T}\textrm{Spin}[\mathbf{R}\boldsymbol{\hat{\sigma}}\mathbf{\hat{n}}]\textrm{d}\Gamma 
\mathbf{G}\delta\mathbf{q}
= \mathbf{0},
\end{equation}
\begin{equation} 
\begin{split}
{\mathbf{R}_\varepsilon}
+\int_{\Gamma}\mathbf{\hat{n}\hat{S}_v}^\textrm{T}\mathbf{U_v}\textrm{d}\Gamma\delta\mathbf{v}
-\int_{\Gamma}\mathbf{\hat{n}\hat{S}_v}^\textrm{T}\mathbf{R}^\textrm{T}\mathbf{U}_\Gamma\textrm{d}\Gamma\mathbf{SG}\delta\mathbf{q}
-\int_{\Gamma}\mathbf{\hat{n}\hat{S}_v}^\textrm{T}\mathbf{R}^\textrm{T}\textrm{Spin}[\mathbf{u}^\textrm{d}_\Gamma]\textrm{d}\Gamma
\mathbf{G}\delta\mathbf{q}
= \mathbf{0}.
\end{split}
\end{equation}
The above linearized equations can are expressed conveniently in the matrix form
\begin{equation}
\left[
\begin{array}{cc}
\mathbf{F} & -\mathbf{A}^\textrm{T}(\mathbf{I}-\mathbf{SG}) - \mathbf{QG} \\
-\mathbf{A} & \mathbf{UG}
\end{array}
\right]
\left[
\begin{array}{c}
\delta\mathbf{v} \\ \delta\mathbf{q}
\end{array}
\right] =
\left[
\begin{array}{c}
-\mathbf{R}_\varepsilon \\
-\mathbf{R}_\sigma
\end{array}
\right],
\end{equation}
where
\begin{equation}
\mathbf{F} =
\int_\Gamma \mathbf{\hat{n} \hat{S}_v}^\textrm{T}\mathbf{\hat{U}_v}\textrm{d}\Gamma,
\end{equation}
\begin{equation}
\mathbf{A} = 
\int_\Gamma \mathbf{\hat{n} \hat{S}_v}^\textrm{T}\mathbf{R}^\textrm{T}\mathbf{U}_\Gamma \textrm{d}\Gamma,
\end{equation}
\begin{equation}
\mathbf{Q} = 
\int_{\Gamma}\mathbf{\hat{n}\hat{S}_v}^\textrm{T}\mathbf{R}^\textrm{T}\textrm{Spin}[\mathbf{u}^\textrm{d}_\Gamma]\textrm{d}\Gamma,
\end{equation}
\begin{equation}
\mathbf{U} = 
\int_{\Gamma}\mathbf{U}_\Gamma^\textrm{T}\textrm{Spin}[\mathbf{R}\boldsymbol{\hat{\sigma}}\mathbf{\hat{n}}]\textrm{d}\Gamma.
\end{equation}


\subsection{Partial Solution}

The size of the problem to be solved can be reduced by static condensation,
taking advantage of the element level approximations of stress fields.
On element level, solution for stress degrees of freedom is expressed by
\begin{equation}
\delta\mathbf{v} = 
-\mathbf{F}^{-1}\mathbf{R}_\varepsilon+
\mathbf{F}^{-1}\left\{
\mathbf{A}^\textrm{T}(\mathbf{I}-\mathbf{SG})
+\mathbf{QG}
\right\}
\delta\mathbf{q}.
\end{equation}
Substituting $\delta\mathbf{v}$ in Eq.~(\ref{eq:res1}) yields
\begin{equation}
\left\{-\mathbf{A}\mathbf{F}^{-1}\left[
\mathbf{A}^\textrm{T}(\mathbf{I}-\mathbf{SG})
+\mathbf{QG}
\right]
-\mathbf{UG}
\right\}
\delta\mathbf{q} = 
-\mathbf{R}_\sigma
-\mathbf{A}\mathbf{F}^{-1}\mathbf{R}_\varepsilon.
\end{equation}
This leads to a modified system of equations
\begin{equation}
\mathbf{\hat{F}}\delta\mathbf{q} = -\mathbf{\hat{R}}_\sigma.
\end{equation}

In the above, the stress degrees of freedom are eliminated conveniently from
the global system of equations, leaving only the displacement degrees of
freedom to be determined. This does not affect the final result in any way and
simply represents a solution technique that reduces the number of equations
that have to be solved simultaneously.  Following the solution of the
displacement degrees of freedom, the stress degrees of freedom can be recovered
on an element by element basis, a process that can be easily parallelized.

It is noted, that independent displacement approximation on elements faces, as
presented here, lead to larger number of unknowns in comparison to classical
FEM. Nevertheless, the matrix sparsity, which has a significant impact on the
efficiency of direct and iterative solvers, is minimised in the presented
approach due to the association of the displacement degrees of freedom with
faces rather than vertices and due to the fact that any given face can only
have two neighbouring elements in either 2D or 3D. For the same reason the
communication between parallel processors is significantly reduced. As shown
in \cite{br}, this approach results in very good scalability.

\section{Numerical examples} \label{sec:examples}

To illustrate features of the CR-HTS approach, two numerical examples are presented.
The first one is an academic problem investigating model's accuracy. The second example further
illustrates model's accuracy for large scale 3D representative model of complex micro structure.

Examples are analysed with the use Portable Extensible Toolkit for Scientific
Computation (PETSc) \cite{petsc} and A Mesh-Oriented datABase (MOAB)
\cite{moab}. Finite elements and CR formulation were implemented in C and C++.

\subsection{Pure bending in 3D}

In the first numerical example the numerical results are compared with an analytical solution.
A beam (0.2x0.2x5.0) with square cross-section is stressed under pure bending. 

In order to get analytical solution, a special case is analysed, i.e. $E=1$ and
$\nu = 0$. A plane curvature is given explicitly as
\begin{equation}
\kappa = \frac{|w{''}|}{(1+{w{'}}^2)^{3/2}}
\end{equation}
and the bending moment according to Bernoulli theory is given by 
\begin{equation}
M =  \textrm{EI} \kappa.
\end{equation}
Strain energy is calculated from 
\begin{equation}
E = \frac{1}{2} \int_{-l/2}^{l/2} \textrm{EI} \kappa^2 \textrm{d}l
\end{equation}

\begin{figure}[!htp]
\centering
\includegraphics[width=1.0\textwidth]{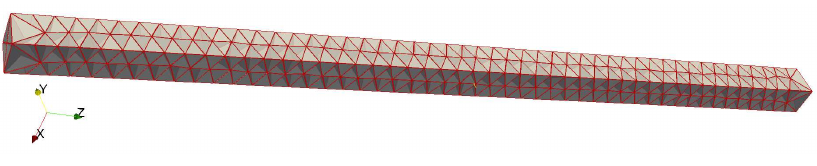}
\caption{ \label{fig:mesh1} Mesh}
\end{figure}
\begin{figure}[!htp] 
\centering
\includegraphics[width=1.0\textwidth]{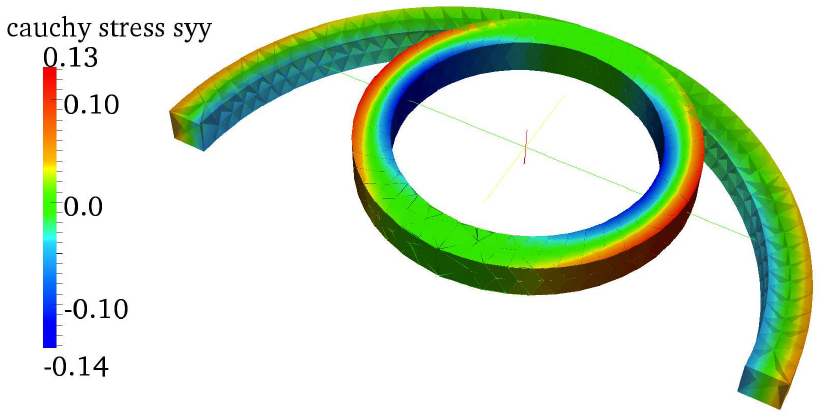}
\caption{\label{fig:cir} Deformation for angles of rotation $\pi$ and $2\pi$ at the end of the bar. The deformation is
plotted using displacements Trefftz approximation functions ($\mathbf{Uv}$) supplemented with rigid body motion.}
\end{figure}
\begin{figure}[!htp]
\centering
\includegraphics[width=1.0\textwidth]{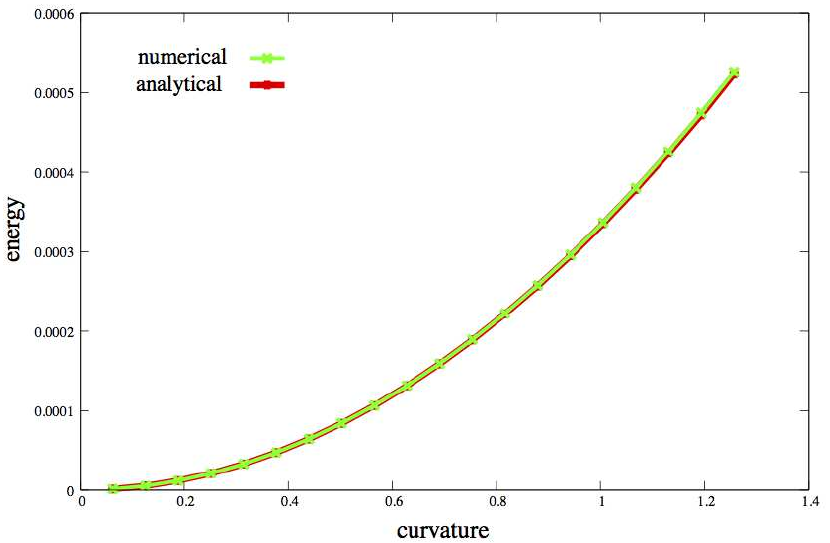}
\caption{Comparison of the analytical and numerical solutions.}
\end{figure}
With the analytical solution at hand, a 3D numerical model is analysed.  A discretized beam consisting 1160
tetrahedrals as shown in Fig.~\ref{fig:mesh1}. A HTS element with quadratic displacements discretization and
3rd order polynomial approximation of stress field within element is used. 
At both ends the kinematic boundary conditions are applied in order to force pure bending
\begin{equation}
u = 0,
\quad v = y \cos \frac{1}{2} L\kappa - y,
\quad w  = \kappa^{-1} \sin \frac{1}{2} L\kappa  - z
\end{equation}
where the length of the beam $L = 5$. Despite the coarse mesh, numerical solution is
in excellent agreement with the analytical one, cf. Fig.~\ref{fig:cir}. For
$\kappa = 2\pi/L$ beam forms a circle, which is in agreement with the analytical
solution. At each time step quadratic convergence is observed.

This test has been reproduced for several bar placements in 3D, in order to verify
the implementation and the formulation of best-fit rotator. We note that the Cauchy stress
field distribution is in an exact much which the analytical solutions for all cases. 

At each increment an updated rotor $\mathbf{R}$ is computed for each element
independently.  At all times the procedure presented in section~\ref{sec:bestrorarir}
required few iterations to achieve the machine precision. Moreover, profiling the
runtimes showed that the time related to the rotor computations was negligible. 

\subsection{Analysis of foam}

\begin{figure}[!htp]
\centering
\includegraphics[width=0.6\textwidth]{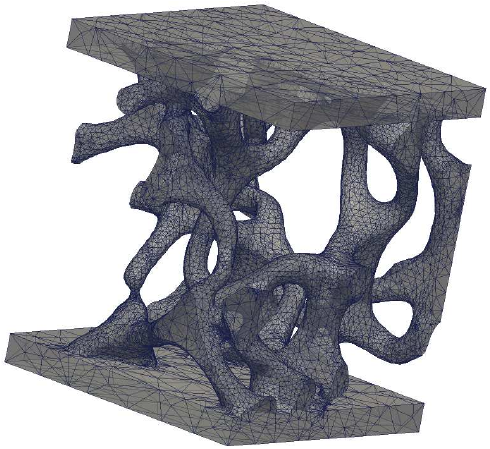}
\caption{Foam surface mesh and geometry}
\end{figure} 
In order to present one of the many potential applications of the Lagrange CR
formulation we investigate deformation of a hypothetical foam. The mesh was
obtained from a high-resolution 3D scan, included as a case study in the Simpleware
software \cite{simpleware}. The mesh comprised 95956 tetrahedrons and 40725
nodes. The constitutive parameters of the model were chosen arbitrarily,
i.e. $E=1e4$ and $\nu=0.25$. At the bottom plate a fixed boundary conditions were
applied. At the top plate tractions in $z$ direction were applied. On the front,
back, left and right side the displacement in the normal direction to the side wall
was blocked.

Displacement and stresses filed were approximated using Hybrid-Trefftz elements.
Two cases were considered with linear and quadratic approximation of
displacements. For the first case, i.e. linear functions approximating
displacements, polynomial Trefftz functions of order two were used to
approximate stresses. For the second case, polynomial Trefftz functions of order
three were used to approximate the stress field. This led to the model with  $\sim 5$
million DoFs and $\sim 10$ million DoFs, for the case one and case two, respectively.
We note, that we not observe shear locking phenomenon for the HTS elements
\cite{jirousek}, despite tetrahedrons being used to approximate the geometry. 

Large rotations led to local buckling in the micro-structure which had an impact
on the convergence of the Newton method. In order to improve efficiency, a set of
nonlinear equations was supplemented with the arc-length control. We introduced
additional unknown $\lambda$, such that 
\begin{equation}
\mathbf{R}_q = \mathbf{\hat{R}}_\sigma - \lambda\mathbf{T}_\lambda,
\end{equation}
and additional control equation
\begin{equation}
\mathbf{R}_\lambda = \mathbf{G}\Delta\mathbf{q}^\textrm{T}\mathbf{G}\Delta\mathbf{q} + \Psi^2\lambda^2 - s^2.
\end{equation}
$\mathbf{G}\Delta\mathbf{q}$ is an increment of the rotation angle $\Delta\vec{\boldsymbol{\Phi}}$, which yields
\begin{equation}
\mathbf{R}_\lambda = \Delta\vec{\boldsymbol{\Phi}}^\textrm{T}\Delta\vec{\boldsymbol{\Phi}} + \Psi^2\lambda^2 - s^2,
\end{equation}
where elements of matrix $\mathbf{G}$ are determined at the beginning of a load
increment. $\Psi$ and $\lambda$ are scaling parameter and load factor
respectively. $s$ is radius of a hyper-sphere. We note that the change of the rotation
angle is a source of nonlinearities and a direct control of the rotation angle
positively affects the stability of the Newton method. Since controlling geometrical
instabilities in a structure like foam is an open problem, the above solution has
been adopted only as a demonstration, while further improvements could well be introduced.

\begin{figure}[!htp] 
\centering
\psfrag{case1}{\small Case 1, linear approximation of displacements}
\psfrag{case2}{\small Case 2, quadratic approximation of displacements}
\psfrag{linear}{\small Linear response}
\psfrag{snap shot}{\small snapshot}
\psfrag{arc}{$s$}
\psfrag{lambda}{$\lambda$}
\includegraphics[width=1\textwidth]{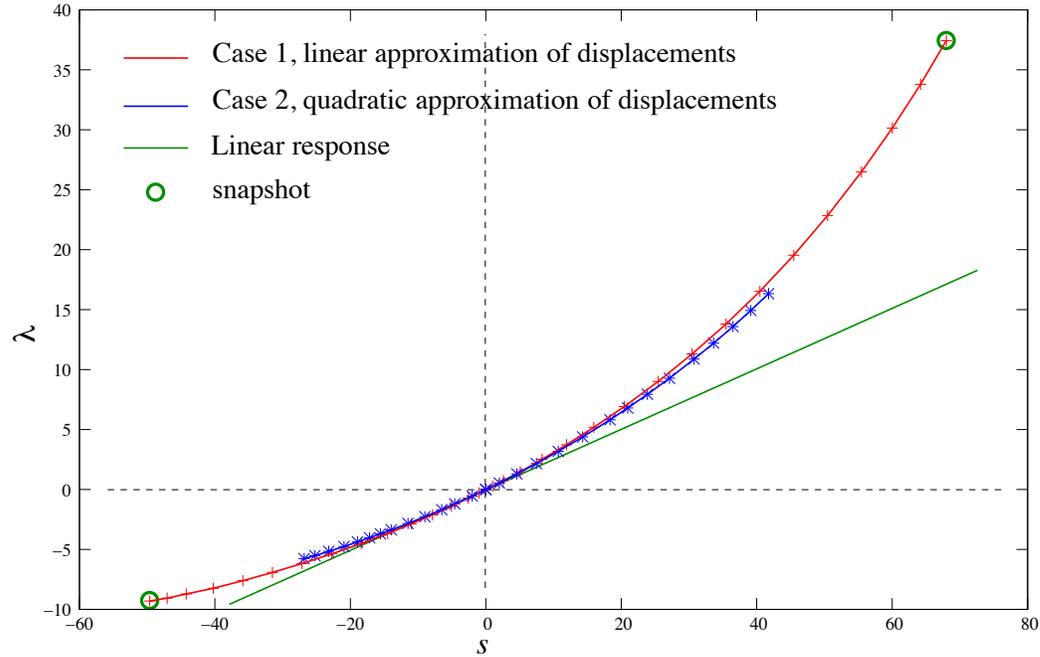} 
\caption{\label{fig:fr0} Equilibrium path}
\end{figure}
\begin{figure}[!htp] 
\centering
\includegraphics[width=0.9\textwidth]{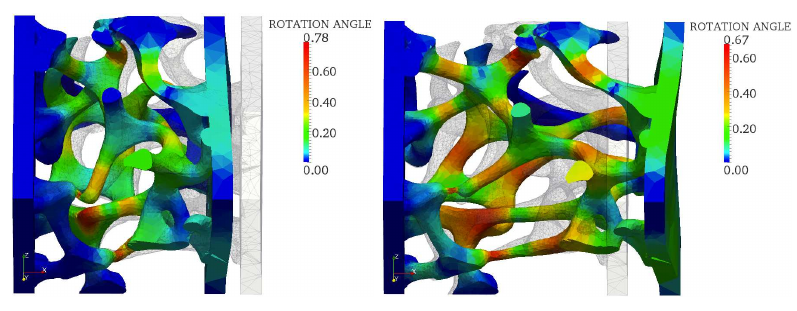} 
\caption{\label{fig:fr1} Rotations, see Fig.~\ref{fig:fr0} where snapshots are taken.}
\end{figure}
\begin{figure}[!htp] 
\centering
\includegraphics[width=0.9\textwidth]{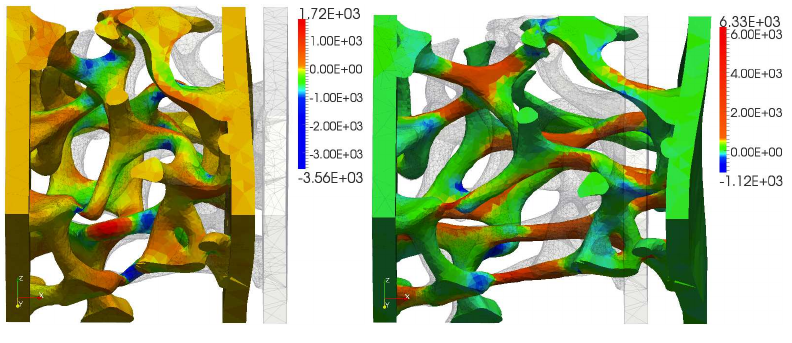}
\caption{\label{fig:fszz1} Stress $\sigma_{xx}$, see Fig.~\ref{fig:fr0} where snapshots are taken.}
\end{figure}

Analysis of foam deformation under compression and tension showed expected foam
response, i.e. overall response softening and stiffening for compression and
tension, respectively. Deformation, rotation angle and stress $\sigma_{xx}$ are
shown in Fig.~\ref{fig:fr1} and Fig.~\ref{fig:fszz1}. We can note that for
compression we observed interpenetration and the model should be extended to take
into account self-contact. For compression, a transition to bending dominated
behavior was observed. For tension, a straightening of foam's branches was observed.
In both cases local buckling was a source of numerical difficulties.

Comparison between linear and quadratic approximation of displacements, i.e.
case 1 and case 2, showed a close much between the numerical responses. Good
quasi-quadratic convergence for case 1 and case 2 was achieved until the final point
of analysis. However poor mesh quality and local instabilities were the source of
problems with the convergence for extreme (non-physical) loads. Those problems
need further investigation and in authors opinion can be overcome by a better local
arc-length control and line searches. 

In figures showing deformation, displacements are plotted using Trefftz
displacement function $\mathbf{u}$, see Eq.~\ref{eq:f7a}, supplemented with
rigid body rotation $\mathbf{u}^r$.  Trefftz displacement field is adjoint to
stresses calculated to fulfill equilibrium, hence displacement vector $\mathbf{u}$ in
Eq.~\ref{eq:f7a} do not carry any information about translations and rotations. By
visual inspection of the figures in this example and in Fig.~\ref{fig:cir} we can verify
an excellent consistency of the calculated rotor with the overall deformation.



\section{Conclusions} \label{sec:concl}

This paper presents a geometrically nonlinear Lagrange co-rotational formulation
applied to Hybrid Stress Trefftz Elements. All terms in the internal force and
tangent matrices are accounted for. 

A new best-fit rotator for continuum elements is developed. The best-fit rotator
functional is constructed by utilizing a boundary integral over a deformable element/body.
Stationary point of the functional lead to the set of three nonlinear equations,
solved by means of the Newton method. Formulation of best-fit rotator, load vector and
tangent matrices was verified by two examples: pure bending test and large deformation of foam. 

Findings from foam analysis can be used to model auxetic and non-auxetic
polymeric foams \cite{foam}, remodelling and growth of tubercular bone \cite{bone}
and fracture of ceramic/quasi-brittle foams.

In authors opinion exciting research opportunities arise form this contribution.
The CR-HTS method can be applied to many problems involving quasi-brittle materials and
large deformations.  We note, that for fracturing materials very often large
rotations are observed for fractured parts. Moreover the current approach can be
applied to many dynamic problems involving large rotations and small deformations.
For models involving in addition large strains, it is possible to use the
CR formulation together with Total Lagrangian or Updated Lagrangian formulations,
where the CR frame could be attached to rotating parts of the model or individual
elements.

Another interesting research opportunity is application of the Lagrange CR formulation
to very large problems, solved with Krylov iterative solver without explicit formation
of the tangent matrix. Matrix free methodology is attractive for memory demanding problems.
We note, that in case of the presented methodology, all integrals can be computed only at the beginning
of an analysis, and the tangent matrix can be calculated with the use of projection matrices only.
As a result, the matrix free methodology can be used efficiently with nonlinear problems and
higher-order elements without the need for integration at each iterative solver step.




\begin{thebibliography}{99}

\topsep=0.0ex
\parsep=0.0ex
\parskip=0.0ex
\itemsep=0.0ex

\bibitem{fraitas} J.A. Teixeira de Freitas, \emph{Formulation of elastostatic
hybrid-Trefftz stress elements}, Comput. Methods Appl. Mech. Engrg. 153 (1998),
pp. 127–151.

\bibitem{Ktefftz} Lukasz Kaczmarczyk, Chris J. Pearce, \emph{A corotational
hybrid-Trefftz stress formulation for modelling cohesive cracks},Comput.
Methods Appl. Mech. Engrg, Volume 198, Issues 15-16, 15 March 2009, Pages
1298-1310.

\bibitem{Rankin} B. Nour-Omid, C.C. Rankin, \emph{Finite rotation analysis and
consistent linearization using projectors}, Comput. Methods Appl. Mech. Engrg.
93 (1991) 353–384.

\bibitem{Crisfield1} M.A. Crisfield, \emph{A consistent corotational formulation for
nonlinear three-dimensional beam element}, Comput. Methods Appl.  Mech. Engrg.
81 (1990) 131–150.

\bibitem{Crisfield2} M.A. Crisfield, G.F. Moita, \emph{A unified co-rotational for
solids}, shells and beams, Int. J. Solids Struc. 33 (1996) 2969–2992.

\bibitem{Felippa} C.A. Felippa, B. Haugen, \emph{A unified formulation of
small-strain corotational finite elements: I. Theory}, Comput. Methods Appl.
Mech. Engrg.  194 (2005) 2285–2335

\bibitem{petsc} Satish Balay and Kris Buschelman and William D. Gropp and
Dinesh Kaushik and Matthew G. Knepley and Lois Curfman McInnes and Barry F.
Smith and Hong Zhang, \emph{{PETSc} {W}eb page}, http://www.mcs.anl.gov/petsc,
2009.

\bibitem{moab} Timothy J. Tautges, \emph{MOAB-SD: integrated structured and
unstructured mesh representation}, Engineering with Computers, 20(3) 2004 286-293

\bibitem{simpleware} \emph{Simpleware Ltd {W}eb page}, http://www.simpleware.com

\bibitem{jirousek} J Jirousek, A Wroblewski, QH Qin, XQ H, \emph{A family of
quadrilateral hybrid-Trefftz p-elements for thick plate analysis}, Comput.
Methods Appl. Mech. Engrg. 127 (1995) 315-344.

\bibitem{foam} Samuel A. McDonald, Ghislain Dedreuil-Monet, Yong Tao Yao,
Andrew Alderson, and Philip J. Withers, \emph{In situ 3D X-ray microtomography
study basic solid state physics comparing auxetic and non-auxetic polymeric
foams under tension}, Phys. Status Solidi B, 1–7 (2010) / DOI
10.1002/pssb.201083975.

\bibitem{bone} L. Kaczmarczyk and C.J. Pearce, \emph{Efficient Numerical
Analysis of Bone Remodelling}, Journal of the Mechanical Behavior of
Biomedical, 2010.

\end{thebibliography}
\end{document}